\documentclass[10pt]{amsart}

\usepackage{amsfonts}
\usepackage{amsxtra}
\usepackage{amscd}
\usepackage{graphicx}
\usepackage{color}

\newtheorem{theorem}{Theoreom}[section]

\newtheorem{proposition}[theorem]{Proposition}
\newtheorem{corollary}[theorem]{Corollary}

\newtheorem{definition}[theorem]{Definition}
\newtheorem{example}[theorem]{Example}

\def\KK{{\mathbb K}}

\def\mm{{\mathbb M}}  \def\PP{{\mathbb P}}
\def\NN{{\mathbb N}}

 \def\tg{{\mathbf{t}}}

\def\Cc{{\mathcal C}} \def\Res{{\mathrm{Res}}}

\def\tg{{\mathbf{t}}}
\def\Xg{{\underline{X}}}\def\DD{{\mathbb D}}

\begin{document}

\title[A matrix-based approach to properness and inversion
  problems]{A matrix-based approach to properness and inversion
  problems for rational surfaces}

\date{\today}

\author{Laurent Bus\'e}
\address{INRIA, 2004 route des Lucioles, B.P. 93, 06902 Sophia-Antipolis cedex, France.}
\email{lbuse@sophia.inria.fr}

\author{Carlos D'Andrea}
\address{Departament d'\`Algebra i Geometria,
Facultat de Matem\`atiques, Universitat de Barcelona,  Gran Via de les Corts Catalanes, 585
08007 Barcelona Spain.} \email{carlos@dandrea.name}

\begin{abstract}
We present a matrix-based approach for deciding if the parameterization of an algebraic space surface is invertible or not, and
  for computing the inverse of the parametrization if it exists.
\end{abstract}

\maketitle


\section{Introduction}
Rational surfaces play an important role in the frame of practical
applications, especially in Computer Aided Geometric Design (see
\cite{HCAGD,PACA} and the references therein). These surfaces can
be parameterized, i.e.~can be seen as the image of a generically
finite rational map
\begin{equation}\label{param}
\begin{array}{crcl}
\phi:&{\KK}^2& \dashrightarrow & {\KK}^{3} \\
&\tg:=(t_1,t_2)&\mapsto&\left(\frac{p_1(\tg)}{q_1(\tg)},\frac{p_{2}(\tg)}{q_2(\tg)},\frac{p_3(\tg)}{q_3(\tg)}\right),
\end{array}
\end{equation}

\noindent where $\KK$ denotes the ambient field, which we assume to
be of characteristic zero. In the sequel we will address the following questions:
\begin{itemize}
\item \emph{Properness problem}: decide if $\phi$, co-restricted to
  its image, is invertible.
\item \emph{Inversion problem}: If $\phi$, co-restricted to its image,
  is invertible, then compute its inverse.
\end{itemize}

Both questions have already been solved theoretically and
algorithmically.
For plane curves, the situation is very well-known: one can relate the
properness and inversion problems to L\"uroth's theorem, and there are
different algorithmic procedures to solve them (see
\cite{HCAGD,PACA,SW,vho}). For space surfaces, there exist some algorithmic
approaches based on $u$-resultants \cite{ChGo} and on Gr\"obner bases
in \cite{Sch}. In
\cite{PDSS} a complete algorithm is given to solve both problems by means of univariate  resultants and GCD computations. 
Our starting point is the resultant matrix-based method presented in
\cite[Chapter $15$]{HCAGD} for inverting a parametrized algebraic
surface, and in \cite[\S 5]{BEM03} where it is used
for computing the inverse image of a point of a parameterization.

\medskip

In this paper, we introduce a general matrix-based approach for dealing with both problems.
We will begin by reviewing the plane curves case; this will clarify our approach and help the reader to understand it.
In section \ref{invv}, we will introduce \emph{inversion matrices} associated to parameterizations.
We will show that their existence implies the properness of a given
parameterization and also that they can be used in order to produce an inverse map in
terms of quotient of determinants of some of its sub-matrices.

\par
The problem of deciding whether a given parameterization is proper or not is not an easy task, hence in general
the construction of such matrices is a non trivial problem.
In \cite{BD} some matrices called \emph{Jacobian matrices} were introduced in order to deal with the properness problem.
This is a first example of a family of inversion matrices; we will briefly recall this construction at the end of section \ref{invv}.
\par
In order to obtain other groups of inversion matrices,
we introduce \emph{implicitization matrices} in section \ref{thm}.
They are essentially square matrices whose determinant is a nontrivial multiple of the
implicit equation of the given
surface (although we do not need to compute the implicit equation in the process).
We will prove that they characterize the properness of a given
parameterization and moreover produce inversion matrices if the parameterization is proper.
Then, we give a list of known algorithms to construct implicitization matrices.
We end the paper with some illustrative examples.

\medskip

We would like to emphasize that the scope of this paper is not to provide a complete algorithm to
solve both properness and inversion problems, but to give general theorems showing that it is possible to test
properness and to get inverse maps of a parameterization by using tools from elimination theory which have been
developed in principle for other purposes: as soon as one has what we call an
inversion matrix (or an implicitization matrix) it is possible to read almost immediately the inverse map from it.
\par
A complete algorithm for computing the inverse map from the parameterization by using our approach requires only the construction of an
implicitization matrix. So, the description of such an algorithm is
only the description of the construction of an implicitization matrix. This is why we do not deal with algorithms in this
paper, they can be found in the literature of matrices from elimination theory.
\par
In section \ref{thm} we give a list of some general constructions yielding implicitization matrices and it turns out
that this
list already covers a lot of cases. For the general situation, no process is yet known (as far as we know),
but one can also try to produce it case by case by computing particular syzygies.

It should be pointed out that our approach is very different from what has been done
in \cite{ChGo,PDSS,Sch} since we provide inversion formulas in terms of quotient of
determinants of sub-matrices of a given matrix instead of producing
expanded rational symbolic expressions. For instance, as a consequence of our results, we will show
that, in the case where resultant matrices give implicitization matrices
(this is always the case for plane curves for instance),
the inversion map can be represented by sub-matrices of a resultant matrix which is
built from the coefficients of the given parameterization; in other words, no symbolic
computations are needed if we stay at this level of representation (note that this
kind of representation in Computer Aided Geometric Design as already been widely used,
for instance to deal with the surface/curve or surface/surface intersection problems).
Moreover, it is interesting to know that certain big polynomials can be
represented as determinants of  certain matrices. In this paper we prove that inversion maps are given
in terms of determinants (as soon as one has an implicitization map).

\medskip

Finally, it should be mentioned that the two main theorems \ref{invthm} and
\ref{impthm} proven in this paper are also valid for rational
parameterizations of hypersurfaces (the proofs work verbatim), that is
to say for rational parameterizations from $\KK^{n-1}$ to $\KK^n$,
with $n\geq 2$, whose closed image is a hypersurface. We chose to
stay in the context of surfaces because of their important applications in Computer
Aided Geometric Design.

\section{Preliminary: inversion of parameterized plane algebraic curves}\label{curves}

Before dealing with the case \eqref{param} of surfaces we first quickly describe
our approach to properness and inversion problems in the case of plane algebraic
curves for the sake of clarity and completeness.

Suppose given a rational plane algebraic curve $\Cc$ parameterized by
$$\phi : \KK^1 \rightarrow \KK^2 : t \mapsto
\left(\frac{p_1(t)}{q_1(t)},\frac{p_2(t)}{q_2(t)}\right)$$
where we can assume, without loss of generality, that $\gcd(p_1,q_1)=\gcd(p_2,q_2)=1$. We moreover assume that $m:=\max(\deg(p_1),\deg(q_1)) \geq 1$ and similarly that
$n:=\max(\deg(p_2),\deg(q_2)) \geq 1$; if it is not the case, then $\Cc$ is a line and the properness and inversion problems are easy.

It is well-known that
$$\Res_{m,n}(p_1(t)-xq_1(t),p_2(t)-yq_2(t))=C(x,y)^{\deg(\phi)}$$
where $C(x,y)$ denotes an implicit equation of the curve $\Cc$ and $\Res_{m,n}(-,-)$ the classical Sylvester resultant. This resultant can be computed as the determinant of the so-called Sylvester matrix $S_{m,n}(p_1(t)-xq_1(t),p_2(t)-yq_2(t))$ which satisfies
the equality

\begin{equation*}
    {}^tS_{m,n}(p_1(t)-xq_1(t),p_2(t)-yq_2(t))
  \left(\begin{array}{c}
t^{m+n-1} \\ t^{m+n-2} \\ \vdots \\ t \\ 1
    \end{array}\right)=
  \left(\begin{array}{c}
t^{n-1}(p_1(t)-xq_1(t)) \\ \vdots \\ t(p_1(t)-xq_1(t)) \\ p_1(t)-xq_1(t) \\ t^{m-1}(p_2(t)-yq_2(t)) \\ \vdots \\ t(p_2(t)-yq_2(t)) \\ p_2(t)-yq_2(t)
    \end{array}\right).
  \end{equation*}

We denote by $\mm$ the sub-matrix of the above Sylvester matrix
 obtained
by erasing its last column. For all
$i=1,\ldots,m+n$, we also denote by $\Delta_i$ the signed determinant
of $\mm$ obtained by erasing the $i^{\mbox{\scriptsize
    th}}$ row. In this way, we have
\begin{equation}\label{uni1}
\Res_{m,n}(p_1(t)-xq_1(t),p_2(t)-yq_2(t))=\sum_{i=1}^{m+n}c_i\Delta_i,
\end{equation}
where the $c_i$'s are the entries of the last column of the Sylvester matrix (and are hence polynomials in $\KK[y]$), i.e.~$p_2(t)-yq_2(t)=\sum_{i=0}^{m+n-1}c_it^{m+n-1-i}$.

\begin{proposition}[{\cite[section 2]{BD}}]\label{invprop} With the above notation we have
  $$\deg(\phi)=1 \Leftrightarrow
\gcd(\Delta_1,\ldots,\Delta_{m+n}) \in \KK \setminus \{ 0 \}.$$
Moreover, if $\deg(\phi)=1$ then for all $i=1,\ldots,m+n-1$
the rational map
  $$\KK^2 \rightarrow \KK^1: (x,y) \mapsto
  \frac{\Delta_i}{\Delta_{i+1}}$$
  is an inversion of $\phi$.
\end{proposition}

It is important to notice that this proposition gives {\it closed universal inversion formulas}, that is to say that the inversion maps can be pre-computed in terms of the coefficient of the polynomials $p_1,q_1,p_2$ and $q_2$. Moreover, if we stay at the level of matrices, no symbolic computations are needed: the inversion formula is just a quotient of two determinants of sub-matrices of the Sylvester matrix which is itself directly built from the coefficients of the polynomials $p_1,q_1,p_2$ and $q_2$. In this way, we can {\it represent} an inversion map by two sub-matrices of a Sylvester matrix; when an inverse image is required we just have to instantiate the corresponding point in these sub-matrices and then compute the quotient of two determinants of {\it numeric} matrices. This is very similar to the fact that the Sylvester matrix {\it represents} the implicit equation of $\Cc$.

Finally, note that we could use in the above process the Bezout matrix instead of the
Sylvester matrix that we chosen for simplicity. For a definition and properties of the
Bezout matrix, see \cite{CLO}.

\begin{example} Consider the following easy example of the unitary circle parameterized by
  $$ \phi: \KK^1 \rightarrow \KK^2 : t \mapsto
  \left(\frac{2t}{1+t^2},\frac{1-t^2}{1+t^2}\right).$$
  The associated Sylvester matrix is
  $$S_{2,2}(2t-x(1+t^2),1-t^2-y(1+t^2))
  =\left( \begin {array}{cccc}
-x&0&-1-y&0\\
2&-x&0&-1-y\\
-x&2&1-y&0\\
0&-x&0&1-y
\end{array} \right)$$
from we extract the matrix
  $$\mm=\left( \begin {array}{ccc}
-x&0&-1-y\\
2&-x&0\\
-x&2&1-y\\\
0&-x&0
\end {array} \right).$$
\\
The $3\times 3$ minors of the matrix $\mm$ are then
$$\Delta_1=2x(y-1), \
  \Delta_2=-2x^2, \
  \Delta_3=-2x(y+1), \
  \Delta_4=2x^2-4(y+1).$$
Their gcd is a constant so we deduce that $\phi$ is proper and we can check that  the all the inversion formulas given in proposition \ref{uni1} are equal, that is to say that
$$\frac{\Delta_1}{\Delta_2}=\frac{\Delta_2}{\Delta_3}=\frac{\Delta_3}{\Delta_4}
\in \mathrm{Frac}(\KK[x,y]/I_\Cc),$$
where $I_\Cc=(x^2+y^2-1)$ is the defining ideal of $\Cc$; for instance:
$$\frac{\Delta_1}{\Delta_2}-\frac{\Delta_2}{\Delta_3}=\frac{2x(y-1)}{-2x^2}-\frac{-2x^2}{-2x(y+1)}=-\frac{x^2+y^2-1}{x(y+1)}=0\in \mathrm{Frac}(\KK[x,y]/I_\Cc).$$
We can also check we got inversion formulas, that is to say that
$$\frac{\psi(\Delta_1)}{\psi(\Delta_2)}=\frac{\psi(\Delta_2)}{\psi(\Delta_3)}=\frac{\psi(\Delta_3)}{\psi(\Delta_4)}=t
\in \KK(t)$$
where $\psi$ is the map $\KK[x,y] \rightarrow \KK(t) : P(x,y) \mapsto P(\frac{p_1(t)}{q_1(t)},\frac{p_2(t)}{q_2(t)}).$
For instance:
$$\frac{\psi(\Delta_1)}{\psi(\Delta_2)}=\left( \frac{1+t^2}{2t}\right)  \left( {1-\frac {1-{t}^{2}}{1+{t}^{2}}}\right)=t.$$
\end{example}

\bigskip

In the following, we extend this approach to the more intricate case of algebraic space surfaces. To this task, we will introduce two distinct notions:

\paragraph{\it Inversion matrices:} these matrices will play the role of the matrix $\mm$ above. As soon as such a matrix exists we will prove that the given parameterization is birational and we will deduce inversion maps similarly to proposition \ref{uni1}.

\paragraph{\it Implicitization matrices:} these matrices will play the role of the Sylvester matrix above. From such a matrix we can characterize the properness of the given parameterization $\phi$ and then produce an inversion matrix which is then used to get inversion maps.

\section{Inversion by means of matrices}\label{invv}

Suppose given a rational surface $S$ parameterized by the map
\eqref{param}.
In this section we develop a matrix-based approach to the inversion problem. To do this we introduce a certain class of matrices associated to parameterizations that we will
call \emph{inversion matrices}. We prove that if such a matrix exists for a given parameterization, then this parameterization is birational and we can derive an inversion map in terms of sub-matrices of this inversion matrix.

From now on, we will turn to projective geometry, so we will assume that our given rational surface $S$ is embedded in $\PP^3$  and parameterized by
 \begin{equation}\label{hparam}
\begin{array}{crcl}
\phi:&{\PP}^2& \dashrightarrow & {\PP}^{3} \\
&\tg:=(t_1:t_2:t_3)&\mapsto&\left(p_1(\tg):p_{2}(\tg):p_3(\tg):p_4(\tg)\right)
\end{array}
\end{equation}
where the $p_i(\tg)$'s are homogeneous polynomials of the same degree.
Also, we will hereafter denote by $\KK(S)$ (resp.~ $\KK(\PP^2)$) the field of rational functions over $S$ (resp.~over $\PP^2$) (see e.g.~\cite[\S I.3.]{Hart}).

\medskip

We recall that a \emph{moving surface} of bi-degree $(m;n)$ is a
bi-homogeneous polynomial in the sets of homogeneous variables
$\tg:=(t_1,t_2,t_3)$ and  $\Xg:=(X_1,X_2,X_3,X_4)$ of degree
$m\geq 1$ and
$n\geq 1$ respectively. Thus a moving surface $M(\tg;\Xg)$ can be written as:
 $$M(\tg;\Xg)=\sum_{|\alpha|=n} A_{\alpha}(\tg)\Xg^\alpha,$$
where $\alpha \in \NN^4$ is a multi-index and $A_{\alpha}(\tg)$ is a
homogeneous polynomial of degree $m$ in $\tg$. We will say that this moving
surface \emph{follows the parameterization $\phi$} if $M$ is identically zero after substituting $X_i$ by $p_i(\tg)$ for $i=1,2,3,4$.

\begin{definition}\label{invmat} A  matrix $\mm$ is called an \emph{inversion
  matrix of $\phi$} if it satisfies the two following conditions~:
\begin{enumerate}
\item[\rm (i)] There exists an integer $m\geq 1$ and a subset
      $V=\{\tg^{\alpha_1},\ldots,\tg^{\alpha_d}\}$ of the monomial
      basis
      $\{ \tg^\alpha, \alpha \in \NN^3 \ |\alpha|=m \}$ such that
\begin{equation*}
{}^t \mm \ \left(\begin{array}{c}
  \tg^{\alpha_1} \\ \tg^{\alpha_2} \\  \vdots \\\tg^{\alpha_d} \end{array}\right) =
\left(\begin{array}{c}
    M_1(\tg;\Xg) \\  \vdots \\ M_{d-1}(\tg;\Xg)
\end{array}\right),
\end{equation*}
\noindent where the polynomials $M_i(\tg;\Xg)$ are moving surfaces following $\phi$ of
bi-degree $(m;d_i)$.
Moreover there exists $(\tg^{\beta_1},\tg^{\beta_2},\tg^{\beta_3})\in
      V^3$ satisfying
      $t_1\tg^{\beta_3}=t_3\tg^{\beta_1}$ and
   $t_2\tg^{\beta_3}=t_3\tg^{\beta_2}$ (and hence $t_1\tg^{\beta_2}=t_2\tg^{\beta_1}$).
\item[\rm (ii)] The rank of $\mm$ over the field $\KK(S)$ is exactly $d-1$.
  \end{enumerate}
\end{definition}

An inversion matrix $\mm$ is thus
    a non-square matrix of size $d\times(d-1)$, where $d \geq 3$.  Note that the condition (ii) they have to
    satisfy can actually be checked explicitly; one way to do this is described in
    \cite[proposition 3.5 and remark 3.6]{BD}.

\medskip

Before stating the main result of
this section we need  to recall carefully what
is the degree of the map $\phi$, map that we now assume to be co-restricted to $S$. Roughly speaking the degree of $\phi$ is the finite number of preimages by
$\phi$ of a sufficiently generic point on $S$. More precisely, the map
$\phi$ induces an injective morphism $\phi^\sharp: \KK(S) \rightarrow
\KK(\PP^2): f \rightarrow f\circ \phi$. Thus $\KK(\PP^2)$ is a
finite extension field of $\KK(S)$ and its degree is, by definition,
the degree of $\phi$; usually this is summarized by the formula
$\deg(\phi):=[\KK(\PP^2):\KK(S)]$.  It follows
that $\phi$ is birational if and only if $\phi^\sharp$ is an
isomorphism, i.e.~$\phi$ has degree 1 (see e.g.~\cite[\S I.4.]{Hart}).

\begin{theorem}\label{invthm} With the above notation, if there exists an inversion matrix $\mm$ of $\phi$ then $\phi$ is birational.
   Moreover,  denoting by
   $\Delta_{\alpha}$ the \emph{signed} minor of $\mm$ obtained by erasing
   the line
   indexed by the monomial $\tg^\alpha$, an inversion of
   $\phi$ is given by the rational map
   $$S\subset \PP^3 \rightarrow \PP^2 : \Xg \mapsto
   (\Delta_{\beta_1}(\Xg):\Delta_{\beta_2}(\Xg):\Delta_{\beta_3}(\Xg)),$$
   where $\beta_1,\beta_2$ and $\beta_3$ are such that
   $t_1t^{\beta_3}=t_3t^{\beta_1}$ and
   $t_2t^{\beta_3}=t_3t^{\beta_2}$.
\end{theorem}
\begin{proof}
  Since we assumed that the rank of
     $\mm$ over $\KK(S)$ is $d-1$ we deduce that the kernel of ${}^t
     \mm$ is generated by the vector
\begin{equation}\label{delta}
\left( \begin{array}{cccccc}
    \Delta_{\alpha_1} & \Delta_{\alpha_2} & \cdots &
     \Delta_{\alpha_d} \end{array}\right).
\end{equation}

One has, for all $i=1,\ldots,d$,
$\phi^\sharp(\Delta_{\alpha_i}(\Xg))=\Delta_{\alpha_i}(\phi(\tg))$
and hence it equals the \emph{signed} $(d-1)\times (d-1)$ minor of
 $\phi^\sharp(\mm)$ obtained by erasing the line indexed by
 $\tg^{\alpha_i}$.  But by definition we have, in $\KK[\tg]$,
    \begin{equation}
{}^t \phi^\sharp(\mm) \ \left(\begin{array}{c}
  \tg^{\alpha_1} \\ \vdots \\ \tg^{\alpha_d} \end{array}\right) =
\left(\begin{array}{c}
    M_1(\tg;\phi(\tg)) \\ \vdots \\ M_{d-1}(\tg;\phi(\tg))
\end{array}\right)=\left(\begin{array}{c}
    0 \\ \vdots \\ 0
\end{array}\right).
\end{equation}
 Consequently both vectors
    $$\left(\begin{array}{ccc} \tg^{\alpha_1} & \cdots &
         \tg^{\alpha_d} \end{array}\right)
     \ \mbox{ and } \left( \begin{array}{ccc} \phi^\sharp(\Delta_{\alpha_1}(\Xg)) & \cdots &
     \phi^\sharp(\Delta_{\alpha_d}(\Xg)) \end{array} \right)$$
 generate the kernel of ${}^t \phi^\sharp(\mm)$ over $\KK(\PP^2)$, and thus equal
     up to the multiplication by an invertible element of
     $\KK(\PP^2)$. In particular we have
     $$\phi^\sharp\left(\frac{\Delta_{\beta_1}(\Xg)}{\Delta_{\beta_2}(\Xg)}\right)=
     \frac{t_1}{t_2} ,\ \phi^\sharp\left(\frac{\Delta_{\beta_1}(\Xg)}
     {\Delta_{\beta_3}(\Xg)}\right)=\frac{t_1}{t_3} \mbox{ and }
   \phi^\sharp\left(\frac{\Delta_{\beta_2}(\Xg)}{\Delta_{\beta_3}(\Xg)}\right)=\frac{t_2}{t_3},$$
 which imply that $\phi^\sharp:\KK(S) \rightarrow \KK(\PP^2)$ is an
isomorphism.

Now let the
     rational map
$\psi:S\subset \PP^3 \rightarrow \PP^2: \Xg \mapsto
(\psi_1(\Xg):\psi_2(\Xg):\psi_3(\Xg))$
be an inverse of $\phi$ and $\psi^\sharp:\KK(\PP^2) \rightarrow
\KK(S)$ its associated field embedding. By definition of $\mm$ we have in $\KK[S]$,
and hence in $\KK(S)$,
\begin{equation}
{}^t \mm \ \left(\begin{array}{c}
  \psi(\Xg)^{\alpha_1} \\ \vdots \\ \psi(\Xg)^{\alpha_d} \end{array}\right) =
\left(\begin{array}{c}
    M_1(\psi(\Xg);\Xg) \\ \vdots \\ M_{d-1}(\psi(\Xg);\Xg)
\end{array}\right).
\end{equation}

\noindent But for all $i$ we have
$$\phi^\sharp(M_i(\psi(\Xg);\Xg))=\phi^\sharp(M_i(\psi^\sharp(\tg);\Xg)))=
M_i(\phi^\sharp\circ\psi^\sharp (\tg),\phi^\sharp(\Xg))=M_i(\tg,\phi(\tg))=0$$
in $\KK(\PP^2)$ and hence $M_i(\psi(\Xg);\Xg)= 0$ in $\KK(S)$.
Therefore we deduce that the vector
 \begin{equation}\label{psi}
\left( \begin{array}{cccccc}
    \psi(\Xg)^{\alpha_1} & \psi(\Xg)^{\alpha_2} & \cdots &
     \psi(\Xg)^{\alpha_d} \end{array}\right)
\end{equation}
 generates the kernel of ${}^t \mm^*$ over $\KK(S)$. On the other hand we
 know that the vector \eqref{delta}
 is also a generator of the kernel of ${}^t \mm^*$ over $\KK(S)$. It follows
that  vectors \eqref{delta}  and \eqref{psi}  must equal up to the multiplication by an
invertible element of $\KK(S)$, and the claimed result is proved.
\end{proof}

In general, the computation of an inversion matrix is not obvious.
We will mainly obtain them from matrices coming from elimination theory that
we will describe in section \ref{thm}. There we will see that in most of the cases
there exist algorithms for constructing inversion matrices.
Moreover, we will also see that we can take advantage of the matrix formulation
we have for producing closed and universal inversion formulas for some classes
of surfaces by using resultant-based matrices.

\par We end by showing how we can deduce
inversion matrices from Jacobian matrices; this was the main subject of \cite{BD}.

\begin{example}[Inversion matrices from jacobian matrices]
In \cite{DE}, it is shown that one can construct a hybrid matrix whose determinant is a nonzero multiple of
the resultant, having all rows except one of Sylvester style. In \cite{BD}, we show that the maximal minors of the Sylvester part of this matrix are
subresultants and that we can solve the inversion problem by using them.

  This parameterization is extracted from \cite[example $P_1$]{PDSS}:

  $$p_1 = \frac{t_1}{t_1+t_2}, \
      p_2 = \frac{t_1^2-t_1+1}{t_2+1}, \
      p_3 = t_1^2+t_2.$$

\noindent By considering $F_1:=(t_1+t_2)X_1-t_1,$
$F_2:=(t_2+1)X_2-(t_1^2-t_1+1)$ and $F_3:=X_3-t_1^2+t_2,$ we get the following
subresultant matrix:
 $$\mm:=\left(\begin{array}{cccccc}
    X_2 - 1  &    1   &   X_2  &    -1 &       0 &     0 \\

X_3     &   0   &   -1   &  -1    &    0    &  0 \\

0    &  X_1 - 1 &   X_1 &      0 &       0 &      0 \\

0  &      0  &    0  &   X_1 - 1 &     X_1 &     0\\

0 &       0 &     0  &     0   &   X_1 - 1 &   X_1
\end{array}\right).$$

\noindent All the maximal minors of this matrix are subresultants.
Using theorem \ref{invthm} we can solve the inverse problem. We obtain:
{\small $$t_1=-\frac{\Delta_{t_1}}{\Delta_1}=\frac{X_1(X_2-1-X_3)}{X_2X_1-1-X_2},
\ \mbox{and} \
t_2=\frac{\Delta_{t_2}}{\Delta_1}=\frac{X_2 X_1 - X_2 - X_1
    + 1 - X_1 X_3 + X_3}{X_2X_1-1-X_2}.$$}
\end{example}

\section{Implicitization matrices}\label{thm}
We keep the notation of section \ref{invv} where we developed a matrix-based approach to the properness and inversion problems. In this section we introduce a new kind of matrices, that we will
call \emph{implicitization matrices}. We will prove that, when it exists, such a matrix characterizes the properness of the parameterization $\phi$. Moreover, if $\phi$ is proper we can extract from it an inversion matrix of $\phi$, as defined in definition \ref{invmat}.

Hereafter we will denote
by $F(X_1,X_2,X_3,X_4) \in \KK[X_1,X_2,X_3,X_4]$ the
implicit equation (which is actually defined up to
multiplication by a non-zero constant in $\KK$) of $S \subset \PP^3$. Recall that it is the homogeneous polynomial of minimal degree such that
$F(p_1(t),p_2(t),p_3(t),p_4(t))\equiv 0$ in $\KK[t_1,t_2,t_3]$; its degree is
the degree of the surface $S$ that we denote by $\deg(S)$.

\begin{definition}\label{impmat} A square matrix $\mm$ is an \emph{implicitization
  matrix of the parameterization $\phi$} if it satisfies the three following conditions~:

\begin{enumerate}
    \item There exists an integer $m\geq 1$ and a subset
      $V=\{\tg^{\alpha_1},\ldots,\tg^{\alpha_d}\}$ of the monomial
      basis
      $\{ \tg^\alpha, \alpha \in \NN^3 \ |\alpha|=m \}$ such that
\begin{equation*}
{}^t \mm \ \left(\begin{array}{c}
  \tg^{\alpha_1} \\ \tg^{\alpha_2} \\  \vdots \\\tg^{\alpha_d} \end{array}\right) =
\left(\begin{array}{c}
    P(\tg;\Xg) \\ M_1(\tg;\Xg) \\  \vdots \\ M_{d-1}(\tg;\Xg)
\end{array}\right),
\end{equation*}

\noindent where polynomials $M_i(\tg;\Xg)$ are moving surfaces following $\phi$ of
bi-degree $(m;d_i)$ and $P(\tg;\Xg)$ is an arbitrary bi-homogeneous
polynomial with positive degree in variables $\Xg$. Moreover there
exists  $(\tg^{\beta_1},\tg^{\beta_2},\tg^{\beta_3})\in
      V^3$ satisfying
      $t_1\tg^{\beta_3}=t_3\tg^{\beta_1}$ and
   $t_2\tg^{\beta_3}=t_3\tg^{\beta_2}$ (and hence $t_1\tg^{\beta_2}=t_2\tg^{\beta_1}$).

\item $\det(\mm)=c.F(X_1,X_2,X_3,X_4)^\delta$ where $\delta\in \NN \setminus \{ 0 \}$ and
  $c\in \KK \setminus \{ 0 \}$.
  \item If $\phi$ is birational, i.e.~$S$ is properly parameterized
    by $\phi$, then $\delta=1$.
  \end{enumerate}
\end{definition}

An implicitization matrix of $\phi$ is hence a square $d\times d$ matrix where $d\geq 3$. Of course the name \emph{implicitization matrix} comes from the
condition (2) in this definition. The following theorem shows that implicitization matrices characterize the properness of the map $\phi$ and moreover yield inversion matrices.

 \begin{theorem}\label{impthm} Let $\mm$ be an implicitization
   matrix of $\phi$ and denote by $\mm^*$ the
    sub-matrix of $\mm$ obtained by erasing its first column. Then the gcd of the maximal minors of $\mm^*$ equals $F^p$ with $p\in
   \NN$.

Moreover the following statements are equivalent:
\begin{enumerate}
\item $\phi$ is birational
\item $p=0$
\item $\mm^*$ is an inversion matrix of $S$.
\end{enumerate}
\end{theorem}

   \begin{proof}

     First the fact that the gcd of the maximal minors of $\mm^*$ is a power
     of $F$ follows immediately from the following equality in $\KK[\Xg]$~:
    \begin{equation}\label{dev}
       \det(\mm)=F^\delta=\sum_{|\alpha|=m}
     c_{\alpha}\Delta_{\alpha} \end{equation}
     where the $c_\alpha$ are the coefficients of the erased column of
     $\mm$, since $F$ is an irreducible and homogeneous polynomial.

     Now suppose that $\phi$ is birational (i.e.~proper). Then we know that
     $\det(\mm)=cF$. Looking at the formula \eqref{dev} we deduce that
     $p=0$ since the $c_\alpha$'s have positive degree (recall that
     $P$ is supposed to have positive degree in variables $\Xg$). This implies
     that the rank of $\mm^*$ over the field $\KK(S)$ is $d-1$ (recall
     that $d$ is the size of the square matrix $\mm$), that is to say
     that $\mm^*$ is an inversion matrix.

     Conversely, assume that $p=0$, that is the gcd of the maximal
     minors of $\mm^*$ is a constant.  Then the rank of
     $\mm^*$ over $\KK(S)$ is $d-1$ and $\mm^*$ is an inversion
     matrix. By theorem \ref{invthm}, $\phi$ is then birational.
\end{proof}

\begin{corollary}
If $\mm$ is an implicitization matrix and $\phi$ is not birational, then $\det(\mm)=F^\delta,$ with $\delta>1.$
\end{corollary}
\begin{proof}
If $\delta=1$ then $p=0$ in theorem \ref{impthm} and we would have that $\phi$ is birational.
\end{proof}

We now give  a (non-exhaustive) list of known  implicitization
matrices; they can be divided into two distinct
groups: the moving surfaces matrices and the resultant
matrices.
\paragraph{\it Moving surfaces matrices.}
All matrices coming from the moving surfaces
method introduced by Sederberg \cite{SC} can be used; they are by
definition implicitization matrices. A lot of recent works have
extended the foundational work of Sederberg. At this time, algorithms to
construct an implicitization matrix of a given parameterization $\phi$ are available
if:
\begin{itemize}
  \item $\phi$ has no base points over $\PP^2$ or $\PP^1\times \PP^1$
    (see \cite{CGZ}),
   \item $\phi$ has l.c.i.~base points (plus some other technical
      conditions) over $\PP^2$ (see \cite{BCD}) or $\PP^1\times \PP^1$
      (see \cite{AHW}),
   \item $\phi$ has no base points on a certain projective toric variety (see
     \cite{DK}).
\end{itemize}
This method of moving surfaces is still under development (see for instance \cite{CCL,JSCC}) and the list above will probably be extended in a near future.

\paragraph{\it Resultant matrices.}
The implicitization problem can be solved using some resultant
computations, as it is illustrated in \cite{BEM03,HCAGD}.
The computation of a resultant often involves the construction of a matrix
which sometimes is an implicitization matrix. At their most typical, resultant matrix give universal formulas for a particular class of parameterizations. They are thus very interesting since they allow the design of a pre-computed inversion formula of a specified class of surfaces. Given a parameterization $\phi$,  resultant matrices which are implicitization matrices of $\phi$ are known if
\begin{itemize}
\item if $\phi$ parameterizes a Steiner surface (see \cite{Jou}, and also \cite{AS})
\item if $\phi$ has no base point on $\PP^1\times \PP^1$ (see \cite{DE} and \cite{DiEm}\footnote{see also \cite{CGZ} where it is explained that the rows of Dixon and Bezoutian matrices are moving planes, and hence that any hybrid combination of them give an implicitization matrix})
\item if $\phi$ has no base point in a certain toric variety (see \cite{KG})
\end{itemize}
In the list above we only mentioned general constructions and it should be pointed out that a lot of constructions adapted to particular cases exist.

%

\section{Some illustrative examples}

Below we exhibit three examples in order to illustrate our matrix-based approach to
the inversion and properness problem.

\subsection{}
This example is taken from \cite{SC}. Consider the following parameterization of a cubic surface
with $6$ base points:
 \[
\begin{aligned}
  p_1 &= t_1^{2} {t_2}+2 t_2^{3}+t_1^{2} {t_3}+4 {t_1} {t_2} {t_3}+4 t_2^{2} {t_3}+3 {t_1}
  t_3^{2}+2 {t_2} t_3^{2}+2 t_3^{3}, \\
  p_2 &= -t_1^{3}-2 {t_1} t_2^{2}-2 t_1^{2} {t_3}-{t_1} {t_2} {t_3}+{t_1} t_3^{2}-2 {t_2}
  t_3^{2}+2 t_3^{3},\\
  p_3 &= -t_1^{3}-2 t_1^{2} {t_2}-3 {t_1} t_2^{2}-3 t_1^{2} {t_3}-3 {t_1} {t_2} {t_3}+2 t_2^{2}
  {t_3}-2 {t_1} t_3^{2}-2 {t_2} t_3^{2},\\
  p_4 &= t_1^{3}+t_1^{2} {t_2}+t_2^{3}+t_1^{2} {t_3}+t_2^{2} {t_3}-{t_1} t_3^{2}-{t_2} t_3^{2}-t_3^{3}.
\end{aligned}
\]
If we pick the following moving planes:
 $$
\begin{array}{l}
M_1:=t_1X_1+t_2X_2+t_3X_3,\\
M_2:=t_1(X_2+X_4)+t_2(2X_2-X_3)+t_3(X_2+2X_4),\\
M_3:=t_1(X_3-X_2)+t_2(-X_1+2X_4)+t_3(X_1-X_2),
\end{array}
$$
we can construct the following matrix indexed by the monomials $t_1,t_2,t_3,$ whose determinant is the implicit equation of the surface:
 \[
\begin{pmatrix}
X_3-X_2&-X_1+2X_4&X_1-X_2\\
-X_2-X_4&X_3-2X_2&-X_2-2X_4\\
X_1&X_2&X_3
\end{pmatrix}.
\]
If we erase the first row of this matrix, Theorem \ref{impthm} tells us that the inverse of this parameterization equals
$(\Delta_{t_1}:\Delta_{t_2}:\Delta_{t_3}),$ which itself equals
 {\small $$
\left(
  \left|\begin{array}{cc}
X_3-2X_2&-X_2-2X_4\\
X_2&X_3
\end{array}\right|:
\left|\begin{array}{cc}
X_2+X_4&-X_2-2X_4\\
-X_1&X_3
\end{array}\right|:
\left|\begin{array}{cc}
-X_2-X_4&X_3-2X_2\\
X_1&X_2
\end{array}\right|
\right).$$}

\subsection{}
The following example appears in \cite{DK}.
Consider the following parameterization:

 \[
\begin{aligned}
p_1 &= t_3^3+t_1t_3^2-t_2t_3^2+t_1t_2t_3-t_1^2t_2-t_1t_2^2,\\
p_2 &= t_3^3+t_1t_3^2-t_2t_3^2-t_1t_2t_3+t_1^2t_2-t_1t_2^2,\\
p_3 &= t_3^3-t_1t_3^2+t_2t_3^2-t_1t_2t_3-t_1^2t_2+t_1t_2^2,\\
p_4 &= t_3^3-t_1t_3^2-t_2t_3^2+t_1t_2t_3-t_1^2t_2+t_1t_2^2.
\end{aligned}
\]
There are two moving planes and one moving quadric of degree one that follow the surface:
 $$
\begin{array}{l}
M_1:=t_1(X_4-X_3)+t_3(X_1-X_2),\\
M_2:=t_1(X_2-X_3+2X_4)+t_2(X_2+X_3)+t_3(-X_2-X_3+2X_4),\\
M_3:=t_1(X_1X_2+X_1X_3)+t_2(X_1X_3-X_1X_4+X_2^2+X_2X_4)+\\ +t_3(-2X_1^2+X_2^2+X_2X_4-X_3X_4+X_4^2).
\end{array}
$$
In this case, the transpose of the following matrix is an implicitization matrix
 {\small \[
\begin{pmatrix}
X_4-X_3&0&X_1-X_2 \\
X_2-X_3+2X_4&X_2+X_3&-X_2-X_3+2X_4 \\
X_1X_2+X_1X_3&X_1X_3-X_1X_4+X_2^2+X_2X_4 & -2X_1^2+X_2^2+X_2X_4-X_3X_4+X_4^2
\end{pmatrix}. \]
}
By erasing the last row, we obtain an inversion
$(\Delta_{t_1}:\Delta_{t_2}:\Delta_{t_3})$ where
$$
\Delta_{t_1}= \left|\begin{array}{cc}
0&X_1-X_2 \\
X_2+X_3&-X_2-X_3+2X_4
\end{array}\right|,$$
$$ \Delta_{t_2}=\left|\begin{array}{cc}
X_3-X_4&X_1-X_2 \\
-X_2+X_3-2X_4&-X_2-X_3+2X_4
\end{array}\right|,$$
$$\Delta_{t_3}=\left|\begin{array}{cc}
X_4-X_3&0\\
X_2-X_3+2X_4&X_2+X_3
\end{array}\right|.$$

\subsection{}
We will compute the inverse of the surface parameterized by
$$\begin{array}{ccccc}
X_1=\frac{t_1^2}{t_1^2+t_2^2+1}& &X_2=\frac{2}{t_1^2+t_2^2+1}&&X_3=\frac{t_1+t_2}{t_1^2+t_2^2+1}.
\end{array}$$
Note that here we are working with ``affine'' variables, i.e. we set $X_4=t_3=1.$
We will use the Dixon formulation for the resultant. In order to do
so, consider the polynomials
$$
F_1=(t_1^2+t_2^2+1)X_1-2t_1,\
F_2=(t_1^2+t_2^2+1)X_2-2t_2,\
F_3=(t_1^2+t_2^2+1)X_3-(t_1^2+t_2^2-1).$$
The Dixon matrix of the resultant of $F_1,F_2,F_3$ is
the $6\times6$ matrix
$$\DD=\left(\begin{array}{cccccc}
X_1& 0& 0& X_1-1& 0& X_1\\
X_2-2& 0& 0& X_2& 0& X_2\\
X_3& -1& -1& X_3& 0& X_3\\
0& -X_2-2X_1+2&2X_1& 0& -2X_3&0\\
-2X_1+2-X_2& 0&-2X_3& 0& X_2&0\\
2X_1& -2X_3&0& 0& X_2&0
\end{array}
\right)$$
\noindent whose columns are indexed by the monomials
$1,t_1,t_2,t_1^2,t_1t_2,t_2^2$ (in that order).
So, it turns out that the transpose of $\DD$ is an implicitization matrix, and we can recover the inverse by deleting any row and considering the ratios
$(\Delta_1:\Delta_{t_1}:\Delta_{t_2}).$ For instance, by deleting the
last row we get
 $$\begin{array}{cccccc}
t_1&=&\frac{\det\left(\begin{array}{ccccc}
X_1 & 0& X_1-1&0& X_1 \\
X_2 - 2 &             0  &        X_2&     0 &    X_2  \\
X_3    &           -1  &        X_3 &    0 &    X_3  \\
 0 &          2 X_1 &      0 &      -2 X_3 &   0 \\
 -X_2 - 2 X_1 + 2&    -2 X_3 &     0 &       X_2 &     0
\end{array}\right)
}{
\det\left(\begin{array}{ccccc}
 0  &          0 &         X_1 -1 &0&    X_1\\
 0   &         0 &         X_2 &    0  &   X_2  \\
  -1  &         -1&          X_3&     0 &    X_3  \\
  -X_2 - 2 X_1 + 2&    2 X_1  &     0 &      -2 X_3 &   0 \\
  0  &         -2 X_3&      0 &       X_2&      0
\end{array}\right)
}, \\
t_2&=&\frac{\det\left(\begin{array}{ccccc}
 X_1&                0 &          X_1 - 1 &     0  &    X_1\\
 X_2 - 2 &             0  &           X_2&        0 &     X_2\\
 X_3     &           -1    &        X_3 &       0 &     X_3\\
 0 &          -X_2 - 2 X_1 + 2&      0  &     -2 X_3 &   0 \\
 -X_2 - 2 X_1 + 2 &         0  &           0  &      X_2 &     0
\end{array}\right)
}{
\det\left(\begin{array}{ccccc}
 0  &          0 &         X_1 -1 &0&    X_1\\
 0   &         0 &         X_2 &    0  &   X_2  \\
  -1  &         -1&          X_3&     0 &    X_3  \\
  -X_2 - 2 X_1 + 2&    2 X_1  &     0 &      -2 X_3 &   0 \\
  0  &         -2 X_3&      0 &       X_2&      0
\end{array}\right)},
\end{array}
$$
that is to say
$$t_1=\frac{-2 X_3 (4 X_3^2  + X_2^2  - 2 X_2)
}{-X_2 (X_2^2  + 4 X_2 X_1 - 2 X_2 - 4 X_3^2 )} \ \mbox{ and } \ t_2=\frac{ 4 X_3(X_2
  + 2 X_1 - 2)
}{X_2^2  + 4 X_2 X_1 - 2 X_2 - 4 X_3^2 }.$$

\end{document}